\newfont{\tt}{cmtt10}

\newcommand{\inv}{ \mbox{\tt inv} }

\documentclass{article}      
\renewcommand{\arraystretch}{1.67}   
\usepackage{graphicx}
\usepackage{amsmath}
\usepackage{amssymb}
\usepackage{centernot}

\begin{document}             

\newcommand{\beginproof}{{\bf Proof: }}

\newcommand{\thisistheendoftheproof}{{\sl {\bf QED.}}}

\newcommand{\ccom}[1]{\mbox{\footnotesize{#1}}}

\setlength{\parindent}{0mm}
 \setlength{\parskip}{4pt}

\newtheorem{theorem}{Theorem}[section]
\newtheorem{lemma}[theorem]{Lemma}
\newtheorem{proposition}[theorem]{Proposition}
\newtheorem{definition}[theorem]{Definition}
\newtheorem{corollary}[theorem]{Corollary}

\title{A Preliminary Report on  Search for Good Examples of Hall's Conjecture.} 

\author{St{\aa}l Aanderaa\\ {\small Department of Mathematics, University of Oslo}
\and Lars Kristiansen\\  {\small Department of Mathematics, University of Oslo}
\\   {\small Department of Informatics, University of Oslo}
\and Hans Kristian Ruud\\  
{\small Department of Mathematics, University of Oslo}}

{\tiny \maketitle   }                

\section{Introduction.}

Consider the equation
\begin{align*}
 x^3 -y^2 = k \tag{*}
\end{align*}
where $x,y\in \mathbb{N}$ and $k\in \mathbb{Z}$.
It is easy to see that (*) has infinitely many solutions   where $k=0$
(let $x= t^2$ and $y=t^3$ where $t$ is a natural number).
It turns out that (*)
only has finitely many solutions  in $x$ and $y$ when $k$ is a given  
integer different from $0$. 
Moreover, it is hard to find solutions of (*) where $k$ is small compared to $x$ and $y$.
Hall's Conjecture states that there exists a constant $C$ such that for any solution of (*) where $k\neq 0$,
we have $C\sqrt{x}< |k|$. For more on Hall's Conjecture, see \cite{jchs} and \cite{mhall}.

Hall's Conjecture is neither proved nor disproved. To shed some light upon the conjecture,
researchers has tried to find solutions of (*) where $0 < |k|< \sqrt{x}$. We will
refer to such  solutions as  {\em good examples of  Hall's Conjecture}, and
we will say that $(x,y,k)$ is a {\em good triplet} when $x,y\in \mathbb{N}$ and 
$0< |x^3-y^2|=|k| < \sqrt{x}$.

This paper is a preliminary report on our search for new good examples of Hall's Conjecture.
We present a new algorithm that will detect all good examples within a given search space. We have implemented
the algorithm, and our executions have
so far found five new good examples.

{\em Acknowledgments:} The authors are grateful for the support they
have received from the Norwegian meta-center for computational science
(Notur).

\section{The Basis for the Algorithm}

In this section we will state some basic definitions and theorems.
In the next section we will explain  our algorithm. 

We use capital Latin letters to denote polynomials, and we use small Latin
letters to denote numbers.

\begin{definition} \label{donaldduck}
The polynomials $B$, $C$, $F$ and $H$ are defined by by
\begin{eqnarray}
B(q, p, x)\quad &=& p^2-q^2 x  \\
 C(q, p, x, y)   &=& p^3- 3pq^2x + 2 q^3y \\
F (q, p, x, y)   &=& 4pC       -3  B^2  \\
H(q, p, x, y)   &=& 9FB  - 8C^2  
\end{eqnarray}
\end{definition}

\begin{lemma}\label{buse}
We have
$$F   = p^4-    (6 p^2x + 8 p qy - 3 q^2x^2)q^2 \label{f}$$
and
$$H  = p^6- (15p^4x+40p^3q y-45p^2q^2x^2 
 + 24pq^3xy+27q^4x^3-32q^4y^2)q^2  \label{h} $$
\end{lemma}

\beginproof
The lemma follows straightforwardly from the definition of the polynomials $F$ and $H$.
\thisistheendoftheproof

\begin{theorem} \label{arnebarnebuse}
We have
\begin{itemize}
\item[(1)] $C \equiv  p^3 \quad\mod q^2$
\item[(2)] $F \equiv  p^4 \quad\mod q^2$
\item[(3)] $H \equiv p^6 \quad \mod q^2$
\item[(4)] $H \equiv -8 C^2\ \quad \mod 9 |F|$.
\item[(5)] $p^4-2 p C + F \equiv  0 \quad \mod q^3$
\item[(6)] $4 p^6-5p^3C + H  \equiv  0 \quad  \mod q^3$.
\end{itemize}
\end{theorem}

\beginproof
Clause  (1) and (4) follow straightforwardly from the definition of the polynomials $B$, $C$ and  $H$.
Clause (2) and (3) hold by Lemma \ref{buse}.
Furthermore, (5) holds since
\[\renewcommand{\arraystretch}{1.6} \begin{array}{lclll}
 p^4-2 p C + F   & \stackrel{\mbox{\tiny (a)}}{=} &   p^4 - 2pC +  4pC-  3B^2   &  &   \\
   & = &   p^4 + 2pC -  3B^2   & &   \\
   & \stackrel{\mbox{\tiny (b)}}{=} &   p^4 + 2pC -  3(p^2-q^2 x)^2   &  &  \\
   & = &   p^4 + 2pC - 3p^4 +6p^2q^2x - 3q^4 x^2  &  &    \\
   & = &   -2p^4 + 2pC  +6p^2q^2x - 3q^4 x^2  &  &    \\
   & \stackrel{\mbox{\tiny (c)}}{=} &   -2p^4 + 2p(p^3- 3pq^2x + 2 q^3y)  +6p^2q^2x - 3q^4 x^2   & & \\
   & = &   -2p^4 +  2p^4  - 6p^2q^2x + 4p q^3y  +6p^2q^2 -3 q^4 x^2   &    \\
  & = &      4p q^3y   - 3q^4 x^2   &    \\
 & = &     (4p y - 3 q x^2)q^3   &    
\end{array} \renewcommand{\arraystretch}{1.0}\]
where the equalities labeled (a), (b) and (c) hold by the definitions of respectively $F$, $B$ and $C$.
The proof of (6) is also straightforward.
\thisistheendoftheproof

The proof of the next theorem is long and involved.
Most of the poof can be found  in the Section \ref{joggetur}.
The reader should note that the $p$ and the $q$ given by the theorem are such that $\frac{p}{q}$ is a rational
approximation to $\sqrt{x}$.

\begin{theorem} \label{bounds}
Let $(x,y,k)$ be a good triplet.
Then, there exists $p,q \in \mathbb{N}$ such that $p < x^{2/3}+1$ and $q < x^{1/6}$
and
\begin{itemize}
\item $0  <  C(q,p,x,y)  <  3 q x^{1/6}+1$ 
\item  $|  F(q,p,x,y) | < 8 q + 1$
 \item  $|  H(q,p,x,y) | <  72 q^4 + 1$.
\end{itemize}
\end{theorem}

The  final theorem in this section is a straightforward consequence of Definition \ref{donaldduck}.

\begin{theorem} \label{oledoledoffen}
We have
\begin{itemize}
\item ${\displaystyle B \;\;  = \;\;  \frac{H+8C^2}{9F}}$
\item ${\displaystyle p \;\;  = \;\;  \frac{F +3B^2}{4C}}$
\item ${\displaystyle x \;\;  = \;\; \frac{p^2-B}{q^2}}$
\item ${\displaystyle y  \;\;  = \;\;  \frac{3pq^2x - p^3 + C}{2q^3}}$.
\end{itemize}
\end{theorem}

\section{The algorithm.}

Our algorithm works by 
examining quadruples $(q, f, c, h)$. 
For every good triplet $(x, y, k)$,  Theorem  \ref{bounds} guarantees
at least one  
quadruple 
$(q, f, c, h)$  such that
\begin{itemize}
 \item[(i)] $q  <  x^{1/6}$
\item [(ii)]$0 < c = C(q,p,x,y)  <  3 q x^{1/6}+1$
\item[(iii)] $f= |  F(q,p,x,y) | < 8 q + 1$
\item[(iv)] $h = |  H(q,p,x,y) | < 72 q^4 + 1$.
\end{itemize}
The algorithm uses the equalities in
 Theorem \ref{oledoledoffen} to compute  good triplets from quadruples.
The algorithm uses the modulo equivalences of Theorem \ref{arnebarnebuse} to find quadruples that may yield a
good triplets.

Choose $x_{\max}$ and set $q_{\max} = x_{\max}^{1/6}$.
The following algorithm finds all good triplets
with $x < x_{\max}$:

\begin{enumerate}
 \item For $q = 2, \ldots, q_{\max}$,
 examine the corresponding values of $f, c, h$ as outlined in the following steps:
\item  Introduce the auxiliary variable $p_0$. 
Examine each $p_0<q$ such that  $p_0$ and $q$ are co-prime.
Let $p$ be such that $p_0 \equiv p \mod q$ (we do not compute $p$). 
By  Theorem \ref{arnebarnebuse},  $f \equiv  p^4 \mod q^2$, so  $f \equiv  p_0^4 \mod q$.
Also, by (iii), $\vert  f \vert <  8 q + 1$. Thus we can describe (and compute)  the set  $S_f$ of possible values  of $f$ by
$$S_f =\{ iq +( p_0 ^4 \mod q) \mid -8\leq i \leq 8 \}\;.$$
\label{algf}
\item For $f \in S_f$, introduce the auxiliary variable $p_1$.
Examine all $p_1$ such that  $p_1^4 \equiv f \mod q^2$ and $0 < p_1 < q^2$.
Define the set 
$$ S = \{ (f \mod q^2) + qi \mid 0\leq i\leq q-1\}\; .$$ 
Then $p_1 \in S $; we find the admissible values for $p_1$ by checking  the elements of  $S$.
There are at most 4 possible values of $p_1$.
Now $p_1 \equiv p \mod q^2$.
\item By  Theorem \ref{arnebarnebuse},  $c \equiv  p^3 \mod q^2$. An upper bound for $c$ is provided by
(ii). Hence, the  possible values of $c$ are

$$ (p_1^3 \mod q^2), \, (p_1^3 \mod q^2) + q^2,\, \ldots, \, (p_1^3 \mod q^2) + q^2 \left\lfloor \frac{q_{\max}}{q}
\right\rfloor \; .$$
\item Introduce the auxiliary variable $p_2$. Examine all $p_2$ such that
 $p_2 \equiv  p_1 \mod q^2$ and  $p_2^4 - 2 p_2 c + f\equiv 0 \mod q^3$ and  $ p_2 < q^3$.
 These $p_2$ satisfies clause (6) in Theorem \ref{arnebarnebuse}.
\item Introduce  $h_2 \equiv 5 c p_2^3 - 4 p_2^6 \mod q^3$.  This implies $h_2 \equiv h \mod q^3$.
\item Introduce  $h_3 \equiv h_2 \mod q^3$ such that $h_3$ satisfies  (4) and (6) in Theorem \ref{arnebarnebuse}.
This means that we must find $h_3$ such that  $h_3 \equiv h_2 \mod q^3$ and $h_3 \equiv (- 8 c^2)\mod (9 f)$.
By the Chinese Remainder Theorem:
\begin{description}
\item if $q \centernot\mid 3 $ then 
\item \quad   $h_3 = h_2 ( 9 \vert f\vert)\inv_{q^3}(9 \vert f \vert) + (-8 c^2) q^3 \inv_{(9\vert f\vert)}(q^3)$
\item\quad $h \equiv h_3 \mod (9 q^3\vert f \vert)$
\item if $q \mid 3 $ then 
\item \quad     $h_3 = h_2 ( \vert f\vert\inv_{q^3}(\vert f \vert) + (-8 c^2) q^3 \inv_{\vert f\vert}(q^3)$
\item\quad $h \equiv h_3 \mod (q^3 \vert f \vert)$
\end{description}
\item Having established the preceding relation for $h$, 
using that  $ 0 < \vert h \vert < 72 q^4$,  
we are now in a position to find the possible values of $h$ (given $q, f, c$). 
\item The final steps of the  algorithm consists of using the equations in Theorem \ref{oledoledoffen} to
compute a good triplet from  $q, f, c, h$ (if such a triplet exists):
\item $b = \frac{h+8 c^2}{9 f}$; if $b$ is not an integer, then we do not have a good triplet.
\item $p = \frac{f+3b^2}{4c}$; if $p$ is not an integer, then we do not have a good triplet.
\item $x = \frac{p^2-b}{q^2}$; if $x$ is not an integer, then we do not have a good triplet.
\item $y = \frac{3 p q^2 x-c}{2 q^3}$; if $y$ is not an integer, then we do not have a good triplet. 
\item $k = x^3 - y^2$; if $\vert k\vert < \sqrt{x}$,  we have a good triplet. 
\end{enumerate}

\subsection*{Validity of the algorithm}
The algorithm works by identifying rational approximations $\frac{p}{q}$ to $\sqrt{x}$, where $x$
is the first component of a good triplet. These approximations are found by
examining tuples $(q, f, c, h)$: for every good triplet $(x, y, k)$
at least one such tuple  $T_x$ exists, with $q < x^{1/6}$
(Theorem  \ref{bounds}).
Now assume that $\frac{p}{q}$ is an approximation to $\sqrt{x}$ where  $x$ yields a  ``good example''.
Also assume that $x < x_{\max}$. 
By the following argument, the algorithm exhausts all
the possibilities for  $(q, f, c, h)$.

For each $q$, we examine all $p_0$ such that $1 \leq p_0 < q$, with $p_0$ and $q$
are co-prime (if $p_0$ and $q$ should have a common factor $r$, then
$\frac{p_0}{q}$  would be equivalent to
$\frac{p_0/r}{q/r}$).

By Theorem \ref{arnebarnebuse} $f \equiv p^4 \mod q^2$. 
Writing $p = r q + p_0$, we get $f \equiv p^4 \equiv 4 r q  p_0^3 + p_0^4 \mod q^2$.
By Theorem \ref{bounds}, $\vert f \vert < 8 q + 1$. 
Thus $f = (p_0^4 \mod q^2) + s q$, with $\vert s \vert \leq 8$ 
(step \ref{algf}
in the preceding description of the algorithm).
 
By Theorem \ref{arnebarnebuse}, $f \equiv p^4 \mod q^2$. 
Let  $S_{p_1}$ be the set of (at most 4) solutions of the quadratic congruence 
$f \equiv p_1^4 \mod q^2$.
Then  $(p \mod q^2) \in S_{p_1}$.

Now, for each $p_1$ in $S_{p_1}$,  let $S_{c,p_1}$ be the set consisting of the values
$$ (p_1^3 \mod q^2), \, (p_1^3 \mod q^2) + q^2,\, \ldots, \, (p_1^3 \mod q^2) + q^2 
\left\lfloor \frac{q_{\max}}{q}\right\rfloor \; .$$
The $c$ corresponding to the tuple $T_x$
 will be in   $S_{c,p_1}$ for one of the  $p_1$ in $S_{p_1}$.

Let $p_2 \equiv p \mod q^3$. 
By clause (5) in Theorem \ref{arnebarnebuse}, we have
 $p_2^4-2 p_2 c + f \equiv  0\quad\quad \mod q^3$.
Define the set $S_{p_2,c, f} = $ the elements of the series $p_1, p_1 + q^2, \ldots $ (up to $q^3$)
that satisfy this clause. Then $p_2 \in S_{p_2,c, f} $, with $c$ and $f$ the corresponding values in
 $T_x$.

By clauses (4) and (6) in Theorem \ref{arnebarnebuse}, we have $h \equiv -8 c^2 \mod (9\, |F|)$ and
 $h \equiv 5 c p^3 - 4 p^6 \mod q^3$.
Introduce the variable $h_3$ = the smallest positive integer that satisfies these clauses; 
and define the set $S_h( c, f)$ as the collection of all values  $ < 72 q^4 + 1$ that satisfies these clauses
(for given $p_0, c, f$).
Then $h \in S_h( c, f)$, and the algorithm will accordingly produce the output of the quantities $b, p, x, y, k$.

\section{Results.}
In order to  investigate the feasibility of the algorithm, the algorithm was implemented in Python
and tested with values of $q_{\max}$ up to 1000.
As results looked promising, the algorithm was reimplemented in C, using the 
Gnu Multi-Precision library to carry out operations with arbitrary-length integers.
This program was run with $q_{\max} = 10 000$ (corresponding to a $x_{\max}$ of $10^{24}$).
This run took 57 processor-days; after 35 days it produced the solution $\sharp 44$ in Table 1.
A subsequent run, with $q_{\max} = 20 000$ (corresponding to a $x_{\max}$ of $64 \times 10^{24}$), took about 441 processor-days
and reproduced a solution earlier found by  Calvo et al. \cite{JCHShomepage}. Currently the program is running on
the Norwegian national computing facilities (Notur). So far, five new good examples 
 have been found. All known good examples are included in Table 1.

\renewcommand{\arraystretch}{1.2}   
\begin{table} \label{taba}
\scriptsize
\begin{tabular}{|c|l|c|c|l|}
\hline
\#&$x$&$r^{\mbox{10)}}$&$\frac{p}{q}$&Comments\\
\hline
1& 2& 1.42&& 1)\\
2&5234 &4.26&$\frac{217}{3}$ & 2), 3) \\
3&8158 &3.76&$\frac{271}{3}$ & 2), 3) \\
4&93844 &1.03&$\frac{919}{3}$ & 2), 3), 9)\\
5&367806 &2.93&$\frac{1213}{2}$ & 2), 3) \\
6&421351 &1.05&$\frac{5193}{8}$ &  2), 3)\\
7&720114 &3.77&$\frac{4243}{5}$ &  2), 3)\\
8&939787 &3.16&$\frac{6786}{7}$ &  2), 3)\\
9&28187351 &4.87&$\frac{90256}{17}$ &  2), 3)\\
10&110781386 &1.23&$\frac{115778}{11}$ &  2), 3)\\
11&154319269 &1.08&$\frac{211183}{17}$ &  2), 3)\\
12&384242766 &1.34&$\frac{176419}{9}$ &  2), 3)\\
13&390620082 &1.33&$\frac{177877}{9}$ &  2), 3)\\
14&3790689201 &2.20&$\frac{430980}{7}$ &  3)\\
15&65589428378 &2.19&$\frac{768313}{3}$ &  4)\\
16&952764389446 &1.15&$\frac{79063817}{81}$ &  4)\\
17&12438517260105 &1.27&$\frac{507863263}{144}$ &  4)\\
18&35495694227489 &1.15&$\frac{1030703950}{173}$ &  4)\\
19&53197086958290 &1.66&$\frac{437617999}{60}$ &  4)\\
20&5853886516781223 &46.60&$\frac{6426898417}{84}$ &  4)\\
21&12813608766102806 &1.30&$\frac{17319173410}{153}$ & 4) \\
22&23415546067124892 &1.46&$\frac{68094518942}{445}$ &  4)\\
23&38115991067861271 &6.50&$\frac{108354409918}{555}$ &  4)\\
24&322001299796379844 &1.04&$\frac{387001980055}{682}$ &  4) 9)\\
25&471477085999389882 &1.38&$\frac{83083668769}{121}$ & 4) \\
26&810574762403977064 &4.66&$\frac{359227383073}{399}$ & 4) \\
27&9870884617163518770 &1.90&$\frac{4524186815567}{1440}$ &  5) \\
28&42532374580189966073 &3.47&$\frac{8386886845023}{1286}$ &  5) \\
29&44648329463517920535 &1.79&$\frac{4603857036361}{689}$ & 5)  \\
30&51698891432429706382 &1.75&$\frac{9318491574937}{1296}$ &  5) \\
31&231411667627225650649 &3.71&$\frac{14649368819024}{963}$ &  5) \\
32&601724682280310364065 &1.88&$\frac{39714194816596}{1619}$ &  5) \\
33&4996798823245299750533 &2.17&$\frac{250164969159375}{3539}$ & 5)  \\
34&5592930378182848874404 &1.38&$\frac{32531865160357}{435}$ &  5) \\
35&14038790674256691230847 &1.27&$\frac{392068197831386}{3309}$ & 5)  \\
36&77148032713960680268604 &10.18&$\frac{633004435512983}{2279}$ &  6) \\
37&180179004295105849668818 &5.65&$\frac{678311009850201}{1598}$ & 6) \\
38&372193377967238474960883 &1.33&$\frac{539307656512279}{884}$ & 5)  \\
39&664947779818324205678136 &16.53&$\frac{3652370552518775}{4479}$ &  5) \\
40&2028871373185892500636155 &1.14&$\frac{11181418791644809}{7850}$ &  6) \\
41&10747835083471081268825856 &1.35&$\frac{42884607802081920}{13081}$ &  7) \\
42&37223900078734215181946587 &1.87&$\frac{46777434586297319}{7667}$ &  5) \\
43&69586951610485633367491417 &1.22&$\frac{72198966044283893}{8655}$ &  8)\\
44&3690445383173227306376634720 &1.51&$\frac{121619570207840431}{2002}$ & 5) \\
45&162921297743817207342396140787 &10.65&$\frac{20237053244197156774}{50137}$ & 8) \\
46&1114592308630995805123571151844 &1.04&$\frac{95524640670266092418}{90481}$ & 9) \\
47&39739590925054773507790363346813 &3.75&$\frac{211515916260522809737}{33553}$ &  8)\\
48&862611143810724763613366116643858 &1.10&$\frac{930889835660831460142}{31695}$ &  8)\\
49&1062521751024771376590062279975859 &1.01&$\frac{1095269810850785984986}{33601}$ & 8) \\
50&6078673043126084065007902175846955 &1.03&$\frac{20224028423712303104623}{259396}$ & 5)  \\
\hline
\end{tabular}
\caption{See Table 2 for comments.}
\end{table}
\renewcommand{\arraystretch}{1.4}   
\normalsize

\begin{table}  \label{tabb}
\begin{tabular}{rl}
1)&This solution is not found by the algorithm presented here.\\
2)&Found by M.Hall \cite{mhall}.\\
3)&Found by Gebel, Peth\" o and Zimmer \cite{gpz}.\\
4)&Found by N.D.Elkies \cite{elkies}.\\
5)&Found by Jim\a'enez Calvo, Herranz and S\a'aez \cite{jchs}.\\
6)&Found by Johan Bosman utilizing the software of \\
  &Jim\a'enez Calvo, Herranz and S\a'aez \cite{jchs}.\\
7)&Found by Jim\a'enez Calvo \cite{JCHShomepage}.\\
8)&Found by the authors of this paper.\\
9)&From the Danilov-Elkies infinite Fermat-Pell family.\\
10)&$r= k/\sqrt{x}$. High values of $r$ indicate that Hall's Conjecture is false.\\
\end{tabular}
\caption{Comments to Table 1.}
\end{table}

\newpage

\section{The Proof of Theorem \ref{bounds}}

\label{joggetur}.

\begin{lemma}\label{nyttaarsaften}
Let $(x,y,k)$ be a good triplet. Then, there exists $\gamma\in  \mathbb{R}$ such that
$$y= x^{3/2}(1+\gamma) \;\; \mbox{ and } \;\;   \frac{-|k|}{2x^{5/2}}\; <\;  \gamma \;  < \; \frac{|k|}{2x^{6/2}}\; .
$$
\end{lemma}

\beginproof
For any $x,y\in  \mathbb{N}$, we have $\gamma\in  \mathbb{R}$  such that $y= x^{3/2}(1+\gamma)$.
For convenience, let $w$ denote $\sqrt{x}$. Then, we have
\begin{align*}
y \;\; = \;\; x^{3/2}(1+\gamma)  \;\; = \;\;  w^3(1+\gamma)\; . \tag{*}
\end{align*}
Furthermore,
we have,
$$\begin{array}{lclll} 
\gamma w^3 &=& w^3(1 + \gamma) - w^3 &\;\;\;\;\;\;\;\;\;\;\;\; & \\
  &=& y - w^3 & & \ccom{(*)} \\
&=&  (y^2-w^6)/(y + w^3) & \\
           &=& (y^2-x^3)/(y + w^3)  & & \ccom{since $w=\sqrt{x}$}\\
           &=& -k/(y+w^3)  & & \ccom{since $x^3-y^2=k$}\\
\end{array}$$
This establishes that $\gamma w^3 = -k/(y+w^3)$, and thus 
\begin{align*}
\gamma \;\; = \;\;  \frac{-k}{(y+w^3)w^3}  \;\; = \;\;  \frac{-k}{yw^3 +w^6} . 
\tag{**}
\end{align*}

Next, note that $y$ cannot equal $w^3$ (if $y = w^3 = x^{3/2}$, then $(x,y,k)$ will
not be a good triplet as  $x^3-y^2=0$). So, we have either $y>w^3$ or $y<w^3$.

Assume that $y>w^3  = x^{3/2}$. Then, since $x^3-y^2=k$, we have $k<0$.
By (**), we have
$$
0 \;\; < \;\; \gamma \;\; = \;\;   \frac{-k}{yw^3 +w^6}  \;\; < \;\;  \frac{-k}{2w^6}\;\; = \;\;  \frac{|k|}{2x^{\frac{6}{2}}}\; .
$$

Assume that $y<w^3  = x^{3/2}$. Then, since $x^3-y^2=k$, we have $k>0$.
Moreover, we have $w^2 < y$ (if $y\leq w^2 = x$, then $(x,y,k)$ will
not be a good triplet as  $x^3-y^2 >  \sqrt{x}$).
Now, by (**), we have
$$
0 \;\; > \;\; \gamma \;\; = \;\;   \frac{-k}{yw^3 +w^6}  \;\; > \;\;   \frac{-k}{w^5 +w^6} \;\; > \;\;   \frac{-k}{2w^5}
\;\; = \;\;  \frac{-|k|}{2x^{5/2}}\; .
$$
\thisistheendoftheproof

\begin{lemma}\label{nyttaarsdag}
Let $(x,y,k)$ be a good triplet. Then, there exist $p,q,Q \in  \mathbb{N}$ and $\delta\in  \mathbb{R}$ such that
(i) $p= q\sqrt{x}(1+ \delta)$, (ii) $x^{1/18} < q<  x^{1/6} < Q$ and
 $$
\mbox{(iii)}\;\;\; \;\;\; \frac{1}{q\sqrt{x}(Q+q)} \;\;\;
 < \;\; \; \left\vert \delta \right\vert \;\; \; <  \;\;\; \frac{1}{q\sqrt{x}Q}\; .\;\;\; \;\;\;\;\;\; \;\;\;
$$
\end{lemma}

\beginproof
First we note that $\sqrt{x}$ is an irrational number when $(x,y,k)$ is a good triplet.
(If $\sqrt{x}$ is a natural number, then $(x,y,k)$ will not be a good triplet as $k=0$.
But  $\sqrt{x}$ is either a natural number or an irrational number. Thus we conclude that  $\sqrt{x}$ is irrational.)

Let $a_0, a_1, a_2,\ldots$ be the coefficients for the simple continued 
fraction for $\sqrt{x}$, that is 
$$\sqrt{x}= \lim_{n\rightarrow \infty}[a_0; a_1,\ldots , a_n]$$
and let $h_i$ and $k_i$ be, respectively, the nominator and the 
denominator of the convergent $[a_0; a_1,\ldots a_i]$, that is
$\frac{h_i}{k_i}= [a_0; a_1,\ldots a_i]$. Then, for any $i\in \mathbb{N}$, we have
$$
\frac{1}{k_i(k_{i}+ k_{i+1})} \;\;  <  \;\; 
\left\vert \frac{h_i}{k_i} - \sqrt{x} \right\vert
                  \;\;  <  \;\; \frac{1}{k_{i} k_{i+1}} 
$$
and $k_i< k_{i+1}$. For more on continued 
fractions, see e.g.\ \cite{Niven}.
Now, pick
the least $j$ such that $k_{j+1}> x^{1/6}$. Let $q=k_j$, let $p=h_j$ and
let  $Q= k_{j+1}$. Then, we have
$$
\frac{1}{q(q+ Q)}  \;\;  <  \;\; \left\vert \frac{p}{q} - \sqrt{x} \right\vert
                  \;\;  <  \;\; \frac{1}{qQ} 
$$
where $q <  x^{1/6} < Q$ (we cannot have $q =  x^{1/6}$ as $x^{1/6}\not\in \mathbb{N}$). 
Next, let $\delta$ be the real number such that $p=q\sqrt{x}(1+\delta)$.
Then, we have
$$
\frac{1}{q(q+ Q)}  \;\;  <  \;\; \left\vert \frac{ q\sqrt{x}(1+\delta)}{q} - \sqrt{x} \right\vert
                  \;\;  <  \;\; \frac{1}{qQ} \; .
$$
Thus
$$
\frac{1}{q(q+ Q)}  \;\;  <  \;\; \left\vert \sqrt{x}\delta \right\vert
                  \;\;  <  \;\; \frac{1}{qQ} \; .
$$
Thus
$$
\frac{1}{q\sqrt{x}(q+ Q)}  \;\;  <  \;\; \left\vert \delta \right\vert
                  \;\;  <  \;\; \frac{1}{q\sqrt{x}Q} \; .
$$
\thisistheendoftheproof

The next proposition corresponds to the first clause of Theorem \ref{bounds}.

\begin{proposition}\label{propc} 
Let $(x,y,k)$ be a good triplet. Then, there exist $p,q \in  \mathbb{N}$  such that
$$ 0 \; \;  <  \;\; C(q,p,x,y) \; \;  <  \;\;  3 q x^{1/6}+1\; .$$
Moreover, $q<  x^{1/6}$ and $p<  x^{2/3}+1$.
\end{proposition}

\beginproof
To improve the readability, we will use $w$ to denote $\sqrt{x}$. First we observe that the two preceding lemmas
yield
$p,q \in  \mathbb{N}$ and $\gamma,\delta \in  \mathbb{R}$
such that
$$\begin{array}{lclll} 
C &=&  p^3- 3pq^2w^2 + 2 q^3y   & &\ccom{def.\ of $C$  } \\
  &=&  p^3- 3pq^2w^2 + 2 q^3[w^3(1+\gamma)] & & \ccom{Lem \ref{nyttaarsaften}} \\
&=&   [qw(1+ \delta)]^3- 3[qw(1+ \delta)]q^2w^2 + 2 q^3[w^3(1+\gamma)]  & & \ccom{Lem \ref{nyttaarsdag}} \\
&=&   q^3w^3(1+ \delta)^3 - 3q^3w^3(1+ \delta) + 2 q^3w^3(1+\gamma) & & \\
&=&   q^3w^3[(1+ \delta)^3 - 3(1+ \delta) + 2(1+\gamma)] &  & \\
&=&   q^3w^3(\delta^3 + 3\delta^2 + 2\gamma) \; . & &  \\
\end{array}$$
Thus, whenever $(x,y,k)$ is a good triplet, we can fix  $p,q \in  \mathbb{N}$ and $\gamma,\delta \in  \mathbb{R}$ such that
$C(q,p,x,y) =q^3w^3(\delta^3 + 3\delta^2 + 2\gamma)$. Moreover, Lemma \ref{nyttaarsdag} states that $q<  x^{1/6}$.
We invite the reader to check that it  follows from  Lemma \ref{nyttaarsdag} that $p<  x^{2/3}+1$.
Next, we will use the bounds given in Lemma \ref{nyttaarsaften} and Lemma \ref{nyttaarsdag}, to prove that
\begin{align*}
0 \;\; < \;\; q^3w^3(\delta^3 + 3\delta^2 + 2\gamma)=C
\tag{*}
\end{align*}
and
\begin{align*}
C=  q^3w^3(\delta^3 + 3\delta^2 + 2\gamma) \;\; < \;\; 3 q w^{1/3}+1 
\tag{**}
\end{align*}
This, will complete the proof of the proposition (since  $w^{1/3}= x^{1/6})$.

We prove (*). It follows from Lemma \ref{nyttaarsdag} (iii) that $|\delta| < \frac{1}{2}$ and, thus,
we have
\begin{align*}
\delta^3 + 3\delta^2 \;\; = \;\;  2\delta^2 + \delta^2(1+\delta) \;\;   > \;\; 0.
\tag{$\dagger$}
\end{align*}
Hence
\renewcommand{\arraystretch}{2.2} 
$$\begin{array}{lclll} 
 C \; = \; q^3w^3(\delta^3 + 3\delta^2 + 2\gamma) &> &  q^3w^3 2\gamma  &\;\;\;\;\;\;\;\;\; &\ccom{($\dagger$)  } \\
  &>& {\displaystyle q^3w^3 2 \frac{-|k|}{2w^5} } & & \ccom{ Lemma \ref{nyttaarsaften}  }  \\
&=& {\displaystyle  q^3 \frac{-|k|}{w^2} }  & & \\
&>& {\displaystyle \frac{-|k|}{w}} & & \ccom{since  $ q^3<  (x^{1/3})^3 = w $ }\\
&>&  -1\; . & & \ccom{since $|k|<w$ } \\
\end{array}$$
\renewcommand{\arraystretch}{1.4} 
(We have $|k|<w$ since $w= \sqrt{x}$ and $(x,y,k)$ is a good triplet.)
Now we have proved $C>-1$,  but $C$ cannot be 0 as we e.g.\ have
$p =  (F +3B^2)/4C$ (see Theorem \ref{oledoledoffen}).
Thus we conclude that $C>0$. This proves (*).

We turn to the proof of (**). By Lemma \ref{nyttaarsdag}, there exists $Q> x^{1/6}= w^{1/3}$ such that
\renewcommand{\arraystretch}{2.2} 
$$\begin{array}{lclll} 
 C & = &  q^3w^3(\delta^3 \; + \; 3\delta^2 \; + \; 2\gamma)  &\;\;\;\;\;\;\;\;\; & \\
  & \leq &  q^3w^3(|\delta|^3 \; + \; 3|\delta|^2 \; + \; 2|\gamma|)  &\;\;\;\;\;\;\;\;\; & \\
  &< & {\displaystyle q^3w^3 \left(|\delta|^3 \; + \;  3|\delta|^2 \; + \; \frac{|k|}{w^6} \right) } & & \ccom{Lemma \ref{nyttaarsaften} } \\
  &< & {\displaystyle q^3w^3\left(\left[ \frac{1}{qwQ} \right]^3 \; + \; 3\left[\frac{1}{qwQ}\right]^2 
\; + \; \frac{|k|}{w^6}\right) } & & \ccom{Lemma \ref{nyttaarsdag}  } \\
&=& {\displaystyle \frac{1}{Q^3} \; + \; \frac{3qw}{Q^2} \; + \; \frac{q^3|k|}{w^3} }& & \\
&<& {\displaystyle \frac{1}{(w^{1/3})^3} \; + \; \frac{3qw}{(w^{1/3})^2} \; + \; \frac{(w^{1/3})^3|k|}{w^3}}  & & \ccom{since $q< w^{1/3}<Q$ }\\
&=& {\displaystyle \frac{1}{w} + 3qw^{1/3} \; + \; \frac{|k|}{w^2}}  & & \\
&=&  {\displaystyle 3qw^{1/3} \; + \; \frac{w + |k|}{w^2} } & & \\
&=&  {\displaystyle 3qw^{1/3} \; + \; \frac{2w}{w^2} } & &  \ccom{since $k<  \sqrt{x}= w$ }\\
&\leq &   3qw^{1/3} \; + \; 1 \; . & & \ccom{since $w\geq 2$ }\\ \\
\end{array}$$
\renewcommand{\arraystretch}{1.4} 
This completes our proof. \thisistheendoftheproof

The next two  proposition correspond to, respectively,  the second and third clause of Theorem \ref{bounds}.
Detailed proofs of these two propositions  will not be included in this preliminary report.

\begin{proposition}\label{propf} 
Let $(x,y,k)$ be a good triplet. Then, there exist $p,q \in  \mathbb{N}$  such that
$$  |F(q,p,x,y)| \; \;  <  \;\;  8 q + 1 .$$
Moreover, $q<  x^{1/6}$ and $p<  x^{2/3}+1$.
\end{proposition}

\beginproof
Use Lemma  \ref{nyttaarsaften} and  Lemma \ref{nyttaarsdag} to prove that there exist $p,q\in  \mathbb{N}$
and $\gamma,\delta \in  \mathbb{R}$ such that
$$F(q,p,x,y) \; = \; q^4 x^2(8 \gamma + 8\gamma\delta + 4\delta^3 + \delta^4)\: .$$
Then, use the bounds given in the two lemmas to prove that the proposition holds. 
\thisistheendoftheproof

\begin{proposition}\label{proph} 
Let $(x,y,k)$ be a good triplet. Then, there exist $p,q \in  \mathbb{N}$  such that
$$  |H(q,p,x,y)| \; \;  <  \;\;   72 q^4 + 1 .$$
Moreover, $q<  x^{1/6}$ and $p<  x^{2/3}+1$.
\end{proposition}

\beginproof
Use Lemma  \ref{nyttaarsaften} and  Lemma \ref{nyttaarsdag} to prove that there exist $p,q\in  \mathbb{N}$
and $\gamma,\delta \in  \mathbb{R}$ such that
$$H(q,p,x,y) \; = \; q^6 x^3 (144 \delta\gamma  -32 \gamma^2 + 40 \delta^3 \gamma +120 \delta^2\gamma
+ 6 \delta^5 + \delta^6)\: .$$
Then, use the bounds given in the two lemmas to prove that the proposition holds. 
\thisistheendoftheproof


\begin{thebibliography}{AaKrRuud}

\bibitem{jchs}
I. Jim\a'enez Calvo, J. Herranz and G. S\a'aez: {\em A new algorithm to search for small $\vert x^3-y^2\vert$ values.} 
Mathematics of Computation,
Volume 78, Number 268, October 2009, pp. 2435-2444

\bibitem{elkies} N.D. Elkies: {\em Rational points near curves and small nonzero $\vert x^2 - y^3 \vert$.} 
Algorithmic Number Theory. Proceedings of ANTS-IV, pp. 33-63. W. Bosma, (ed).
Springer, 2000; Lecture Notes in Comput. Sci. 1838.

\bibitem{mhall}
M. Hall: {\em The Diophantine equation  $\vert x^3-y^2\vert$.}
Computers in Number Theory.  A. Atkin and B. Birch (eds.), Academic Press, 1971, pp. 173-198.
\bibitem{gpz} J. Gebel, A. Peth\"{o} and H. G. Zimmer: {\em On Mordell’s equation.}  Compositio Math. 
{\bf 110 } (1998), pp. 335-367


 \bibitem{Niven}
Niven, Zuckerman and Montgomery:  
 {\em An Introduction to the Theory of Numbers.}
John Wiley \& Sons, Inc. 



\bibitem{JCHShomepage} 
 {\verb!http://ijcalvo.galeon.com/hall.htm!}
\\\end{thebibliography}
\end{document}